# Probability density function of the Cartesian x-coordinate of the random point inside the hypersphere

Argyn Kuketayev

4/26/2013

## Abstract

Consider randomly picked points inside the n-dimensional unit hypersphere centered at the origin of the Cartesian coordinate system. The Cartesian coordinates of the points are random variables, which form an n-dimensional vector for each point. Observing only the x-coordinate I obtained its probability density function (PDF) $f_n(x)$. I show that it is related to the Gaussian distribution: in limit the function $g_n(z) = f_n(z/\sqrt{n+2})/\sqrt{n+2}$ converges to the PDF of the standard normal distribution.

## Hyperspheres

The n-dimensional unit hypersphere is defined as a set of points with Cartesian coordinates $(x_1, x_2, \ldots, x_n)$ such that $\sum_{i=1}^{n} x_i^2 \leq 1$. Pick a random point inside the hypersphere and observe its $x_1$ coordinate, which is a random variable. Random point picking inside the hypersphere as well as on its surface is an interesting practical problem in itself and has been studied before[1]. The purpose of this paper is to study the probability distribution of the Cartesian $x_1$ coordinate of the point that was already randomly picked inside the hypersphere somehow.

The volume of n-dimensional unit hypersphere can be computed as follows:

$$V_n = \int_{-1}^{1} V_{n-1} \cdot \left(\sqrt{1-x_1^2}\right)^{n-1} dx_1 \qquad (1)$$

It is easy to understand this equation in the following example. Let us compute the volume of an ordinary 3-dimensional sphere. The sphere can be dissected into a stack of thin cylinders of varying diameters from 0 to 1 then back to 0. The volume of the sphere is a sum of volumes of the cylinders as follows:

$$V_3 = \int_{-1}^{1} \pi r_x^2 dx = \int_{-1}^{1} \pi(1-x^2) dx = \frac{4}{3}\pi,$$

where each cylinder has a radius equal to $r_x = \sqrt{1-x^2}$, and the height equal to $dx$. As expected we got the familiar equation for a volume of a sphere.

Let us consider one more example: a 2-dimensional hypersphere, i.e. an ordinary disk. We are more accustomed to refer to the volume of a two-dimensional object as the *area*. A disk can be seen as a stack of thin rectangles of varying lengths from 0 to 1 then back to 0. The area of a disk is the sum of areas of these rectangles:

$$V_2 = \int_{-1}^{1} 2r_x dx = \int_{-1}^{1} 2\sqrt{1-x^2} dx = \pi,$$

where the width of each rectangle is $2r_x = 2\sqrt{1-x^2}$, and the heights are $dx$.

---

[1] See [Barthe] and [Marsaglia].

Note, that in n-dimensional space (n-1)-dimensional hypersphere will project into a hyperdisk, which has volume equal to $V_{n-1} \cdot \left(\sqrt{1-x_1^2}\right)^{n-1}$, where $V_n$ is the volume of the unit hypersphere in n-dimensional space. In our two examples the three-dimensional sphere and the two-dimensional disk were representing hyperspheres, while the cylinders and rectangles were representing the hyperdisks, i.e. (n-1)-dimensional hyperspheres projected into n-dimensional space. The equation for $V_n$ is well-known[2]:

$$V_n = \frac{\pi^{\frac{n}{2}}}{\Gamma\left(\frac{n}{2}+1\right)} \qquad (2)$$

## Probability Densities

Let us observe the Cartesian $x_1$-coordinate of a point that was randomly picked inside the hypersphere. The value of $x_1$ is a random variable. In the geometrical framework, which was described earlier, this point will belong to one of the thin hypercylinders. Hence, the probability density function (PDF) of $x_1$ can be obtained from the ratio of the volume of the hypercylinder, to which $x_1$ belongs, and the total volume of the hypersphere:

$$f_n(x_1) dx_1 = \frac{V_{n-1} \cdot \left(\sqrt{1-x_1^2}\right)^{n-1} dx_1}{V_n}$$

Further, by plugging the equations for $V_n$ we get the closed-form expression for this PDF:

$$f_n(x) = \begin{cases} \frac{\Gamma\left(\frac{n}{2}+1\right)}{\Gamma\left(\frac{n+1}{2}\right)} \frac{\left(\sqrt{1-x^2}\right)^{n-1}}{\sqrt{\pi}} & , \forall\, x \in [-1,1] \\ 0 & , \forall\, x \notin [-1,1] \end{cases} \qquad (3)$$

Note, that a zero-dimensional hypersphere is a point. Its volume is zero, the same as a volume of zero-dimensional hypercube, which is also a point. A one-dimensional hypersphere in one-dimensional space is simply an interval $[-1,1]$ and its PDF is given by an ordinary uniform distribution: $f_1(x) = \frac{1}{2}$.

In a slightly more interesting case of a two-dimensional hypersphere on a two-dimensional plane, i.e. a disk, the PDF (shown on Figure 1) is given by the equation:

---

[2] E.g. see [Пеньков] p.259, Equation 1.16

$$f_2(x) = \frac{2\sqrt{1-x^2}}{\pi}, \forall\, x \in [-1,1]$$

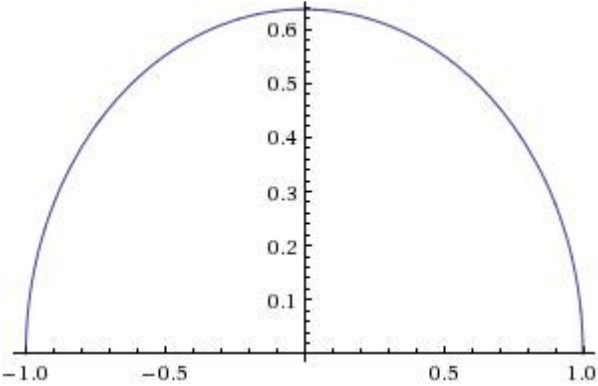

Figure 1. PDF of Cartesian x-coordinate of a random point inside a disk: $f_2(x) = \frac{2\sqrt{1-x^2}}{\pi}$.

## Convergence to Gaussian distribution

The PDFs from Equation 3 for dimensions from 1 to 30 are shown on Figure 2 and Figure 3. It is easy to notice that $f_n(x)$ seems to concentrate around zero as $n$ increases, and that its shape takes a form of a bell-shaped distribution. Indeed I show below that $f_\infty(x)$ collapses to a point and that the PDF of a random variable $z = \sqrt{n+2}\, x$ converges to the standard normal distribution.

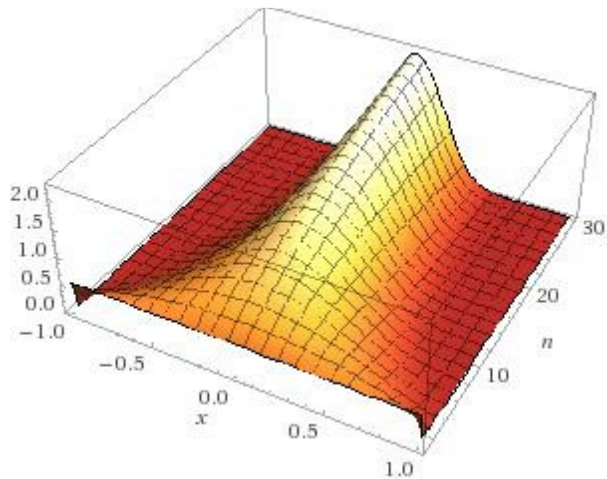
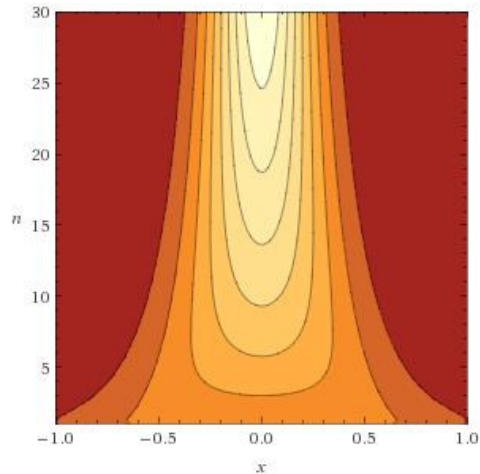

Figure 2. 3D plot of $f_n(x)$, PDF of Cartesian x-coordinate, for n=1 to 30, see Equation 3.

Figure 3 Contour plot of $f_n(x)$, PDF of Cartesian x-coordinate, for n=1 to 30

Consider the PDF of a random variable $z$:

$$g_n(z) = \begin{cases} \dfrac{f_n\left(\frac{z}{\sqrt{n+2}}\right)}{\sqrt{n+2}} & ,\forall z \in \left[-\sqrt{n+2}, \sqrt{n+2}\right] \\ 0 & ,\forall z \notin \left[-\sqrt{n+2}, \sqrt{n+2}\right] \end{cases} \qquad (4)$$

Its characteristics function[3] can be obtained as follows:

$$\varphi_{n,Z}(t) \equiv E[e^{itZ}] = \int_{-\sqrt{n+2}}^{\sqrt{n+2}} g_n(z) e^{itz} dz$$

$$= \int_{-\sqrt{n+2}}^{\sqrt{n+2}} \frac{\Gamma\left(\frac{n}{2}+1\right)}{\Gamma\left(\frac{n+1}{2}\right)} \frac{\left(1-\left(\frac{z}{\sqrt{n+2}}\right)^2\right)^{\frac{n-1}{2}}}{\sqrt{\pi}\sqrt{n+2}} e^{itz} dz$$

$$= \int_{-\sqrt{n+2}}^{\sqrt{n+2}} \frac{\Gamma\left(\frac{n}{2}+1\right)}{\Gamma\left(\frac{n+1}{2}\right)} \frac{(n+2-z^2)^{\frac{n-1}{2}}}{\sqrt{\pi}\sqrt{n+2}\left(\sqrt{n+2}\right)^{n-1}} e^{itz} dz$$

$$\varphi_{n,Z}(t) = \frac{\Gamma\left(\frac{n}{2}+1\right)}{\Gamma\left(\frac{n+1}{2}\right)(n+2)^{\frac{n}{2}}\sqrt{\pi}} \int_{-\sqrt{n+2}}^{\sqrt{n+2}} (n+2-z^2)^{\frac{n-1}{2}} e^{itz} dz \qquad (5)$$

The integral in the last equation can be expressed through the Bessel function of the first kind $J_\nu$ and gamma functions as follows[4]:

$$\int_{-\sqrt{n+2}}^{\sqrt{n+2}} (n+2-z^2)^{\frac{n-1}{2}} e^{itz} dz = \sqrt{\pi}\,\Gamma\left(\frac{n+1}{2}\right)\left(\frac{2\sqrt{n+2}}{t}\right)^{\frac{n}{2}} J_{\frac{n}{2}}\left(\sqrt{n+2}\,t\right) \qquad (6)$$

Plug Equation 5 into the Equation 6:

$$\varphi_{n,Z}(t) = \frac{\Gamma\left(\frac{n}{2}+1\right)}{(n+2)^{\frac{n}{2}}} \left(\frac{2\sqrt{n+2}}{t}\right)^{\frac{n}{2}} J_{\frac{n}{2}}\left(\sqrt{n+2}\,t\right)$$

$$= \Gamma\left(\frac{n}{2}+1\right)\left(\frac{2}{\sqrt{n+2}\,t}\right)^{\frac{n}{2}} J_{\frac{n}{2}}\left(\sqrt{n+2}\,t\right)$$

Finally, we get a simple expression for the characteristic function using the relationship between Bessel function and hypergeometric function[5]:

$$\varphi_{n,Z}(t) = {}_0F_1\left(;\frac{n}{2}+1;-(n+2)\frac{t^2}{4}\right) \qquad (7)$$

---

[3] [Billingsley], p.342

[4] It can be obtained by substituting $a = \sqrt{n+2}$, $\beta = \frac{n+1}{2}$ into an integral $\int_{-a}^{a} (a^2-x^2)^{\beta-1} e^{i\lambda x} dx = \sqrt{\pi}\,\Gamma(\beta)\left(\frac{2a}{\lambda}\right)^{\beta-\frac{1}{2}} J_{\beta-\frac{1}{2}}(a\lambda)$, see Equation 2.3.5 (3) in [Prudnikov].

[5] [Luke], p.332, Section 9.2, Equation (1): $\qquad J_\nu(z) = \frac{(z/2)^\nu}{\Gamma(\nu+1)}\,{}_0F_1(;1+\nu;-z^2/4)$

Using the series representation of $_0F_1$ function[6] we get the Taylor expansion of our characteristic function:

$$\varphi_{n,Z}(t) = \sum_{k=0}^{\infty} \frac{(-(2+n)\,t)^k}{4^k k! \left(\frac{2+n}{2}\right)_k}$$

Let us re-arrange the terms as follows:

$$\varphi_{n,Z}(t) = \sum_{k=0}^{\infty} \frac{(-(2+n)\,t^2)^k \Gamma\left(\frac{n}{2}+1\right)}{4^k k! \, \Gamma\left(\frac{n}{2}+1+k\right)}$$

$$= \sum_{k=0}^{\infty} \frac{(-t^2)^k}{k!\, 2^k} \frac{(2+n)^k \Gamma\left(\frac{n}{2}+1\right)}{2^k \Gamma\left(\frac{n}{2}+1+k\right)}$$

$$= \sum_{k=0}^{\infty} \frac{(-t^2)^k}{k!\, 2^k} \left[ \frac{\left(\frac{n}{2}+1\right)^k \Gamma\left(\frac{n}{2}+1\right)}{\Gamma\left(\frac{n}{2}+1+k\right)} \right]$$

Let us show that each term in the square brackets converge to one as *n* goes to infinity. Using the relation[7]:

$$\Gamma(z+n) = \Gamma(z) z^k \prod_{k=0}^{n-1}\left(1+\frac{k}{z}\right)$$

evaluate the limit:

$$\lim_{n\to\infty} \frac{\left(\frac{n}{2}+1\right)^k \Gamma\left(\frac{n}{2}+1\right)}{\Gamma\left(\frac{n}{2}+1+k\right)} = \lim_{n\to\infty} \frac{\left(\frac{n}{2}+1\right)^k \Gamma\left(\frac{n}{2}+1\right)}{\left(\frac{n}{2}+1\right)^k \Gamma\left(\frac{n}{2}+1\right) \prod_{m=0}^{k-1}\left(1+\frac{m}{\frac{n}{2}+1}\right)} = 1$$

Hence, the equation for characteristic function:

$$\lim_{n\to\infty} \varphi_{n,Z}(t) = \sum_{k=0}^{\infty} \frac{(-t^2)^k}{2^k k!} \qquad (8)$$

---

[6] [Luke], Section 5.2.1, Equation 2: $_0F_1(;b;z) = \sum_{k=0}^{\infty} \frac{z^k}{(b)_k k!}$, where $(b)_k = \frac{\Gamma(b+k)}{\Gamma(b)}$.

7 [Luke], Section 1.1, Equation 2: [Luke]$\Gamma(z+n) = \Gamma(z) \prod_{k=0}^{n-1}(z+k)$

Recall, that the right hand side of Equation 8 is the series representation of the characteristic function of the standard normal distribution[8]: $\varphi(t) = e^{-t^2/2}$.

Therefore, the PDF itself in Equation 4 converges to the PDF of the standard normal distribution:

$$g_n(z) \xrightarrow[n\to\infty]{} \mathcal{N}(0,1) \qquad (9)$$

We also write *informally*[9]:

$$f_n(x) \xrightarrow[n\to\infty]{} \mathcal{N}(0, \frac{1}{n+2}) \qquad (10)$$

The Equation 10 indicates that $\lim_{n\to\infty} f_n(x)$ could probably be a representation of the Dirac's delta-function[10] $\delta(x)$.

## Conclusion

I obtained the PDF of the Cartesian $x$-cooridnate of a random point picked inside the unit n-dimensional hypersphere, and showed that its shape gets closer and closer to the shape of the Gaussian distribution as the dimension of the space increases. This relationship to the normal distribution, as captured in Equation 9, is similar to the central limit theorem[11], and can be interpreted as follows. Then higher is the dimension of the space then less likely it is to observe a Cartesian coordinate of a random point inside the hypersphere equal to the radius of this sphere. In other words, in the infinite-dimensional space all the mass of the hypersphere is concentrated in its center.

Note, that the volume of the n-dimensional hypersphere, given by Equation 2, rapidly shrinks relative to its smallest enclosing hypercube's volume as n increases. This makes the distribution of x-coordinates of the points inside the hypersphere collapse to zero. In other words, a hypersphere from the infinitely dimensional world is projected into one dimensional world as a single point.

---

[8] [Billingsley], Section 26.8
[9] The convergence is not well defined when the right hand side itself is changing, hence the Equation 9 is a proper way of expressing this convergence.
[10] [Li]
[11] [Billingsley], Section 27